\documentclass[11pt]{article}
\usepackage{latexsym,a4}

\usepackage{epsf}
\usepackage{color}
\usepackage{epsfig}

\usepackage{graphicx}
\usepackage{xcolor}
\usepackage{amsmath}
\usepackage{latexsym}
\usepackage{amssymb}
\usepackage{amsfonts}

\usepackage{mathdots}
\usepackage{array}

\usepackage{a4}
\usepackage{longtable}

\usepackage{graphicx}
\usepackage{mathdots}
\usepackage{array}

\title{Parallel Computation of functions of matrices and their action on vectors}

\author{
Sergio Blanes\footnote{Instituto Universitario de Matem\'atica Multidisciplinar,
  Universitat Polit\`ecnica de Val\`{e}ncia, E-46022  Valencia, Spain.
 {\tt serblaza@imm.upv.es}}
}

\makeatletter
\def\Ddots{\mathinner{\mkern1mu\raise\p@
\vbox{\kern7\p@\hbox{.}}\mkern2mu
\raise3\p@\hbox{.}\mkern2mu\raise5\p@\hbox{.}\mkern1mu}}
\makeatother
\usepackage{booktabs}

\def\e{e}

\def\C{\mathbb{C}}

\def\beq{\begin{equation}}
\def\eeq{\end{equation}}

\begin{document}

\maketitle

\begin{abstract}
We present a novel class of methods to compute functions of matrices or their action on vectors that are suitable for parallel programming. Solving appropriate simple linear systems of equations in parallel (or computing the inverse of several matrices) and with a proper linear combination of the results, allows us to obtain new high order approximations to the desired functions of matrices. An error analysis to obtain forward and backward error bounds is presented. The coefficients of each method, which depends on the number of processors, can be adjusted to improve the accuracy, the stability or to reduce round off errors of the methods. We illustrate this procedure by explicitly constructing some methods which are then tested on several numerical examples.
\end{abstract}




\section{Introduction} \label{Introduction}

We present novel algorithms to compute functions of matrices and their action on vectors which can be computed in parallel. We believe that new algorithms that are designed to be computed in parallel will become more frequent and useful in the near future. For instance, about 25 five years ago, in one of the most relevant books in numerical methods \cite{press97nr} it is written:
"In recent years we Numerical Recipes authors have increasingly become
convinced that a certain revolution, cryptically denoted by the words “parallel
programming,” is about to burst forth from its gestation and adolescence in the
community of supercomputer users, and become the mainstream methodology for
all computing...

Scientists and engineers have the advantage that techniques for parallel computation
in their disciplines have already been developed. With multiprocessor
workstations right around the corner, we think that now is the right time for scientists
and engineers who use computers to start thinking parallel." See also e.g. \cite{auckenthaler10mea,blanes21npi,chehab16pmf,gallopoulos16pim,gallopoulos89otp,gander13pap}.


The performance of a parallel algorithm will depend on the particular code in which it is written (Fortran, C++, Matlab, Python, Julia, etc.), the compiler used, the number of processors available, how efficiently they communicate, etc., and this will change in the future. For this reason, we will mainly focus on the structure of the algorithm rather than in its particular implementation for solving a given problem. Then, we will analyse the accuracy and stability of the methods with an error bound analysis and, in some cases, we will show the performance of the new algorithms in the two extreme cases, i.e. in the worst scenario when it is evaluated as a sequential method versus the ideal one in which the cost of the whole algorithm can be taken as the cost for the computations in one single processor and neglecting the cost in the communication between processors. These results will be compared with the results obtained with sequential algorithms from the literature to solving the same problem. These results will illustrate the benefits one can achieve when the parallel algorithms are used under different conditions.


The efficient computation of an important number of functions of matrices of moderate size is of great interest in many different fields 
\cite{almohy09ans,almohy15naf,almohy22apa,arioli96tpm,bader17aia,bader19ctm,bader22aea,golub93mc,higham05tsa,higham09tsa,higham08fom,higham10cma,moler03ndw,najfeld95dot,sastre18eeo,sastre19btc,sastre14aae,seydaoglu21ctm,sidje98eas}. Frequently, it suffices to compute their action on a vector \cite{almohy11ctaf,gallopoulos89otp,higham10cma,hochbruck97oks,hochbruck10ei,niesen12aak} allowing to solve problems of large dimensions, or using an appropriate filtering technique the previous methods can also be used to compute functions of large sparse matrices       \cite{wu21hpc}. For example, exponential integrators have shown to be highly efficient numerical schemes to solve linear systems of differential equations \cite{auckenthaler10mea,blanes09tme,blanes17hoc,gallopoulos89otp,gallopoulos92eso,hochbruck10ei,iserles00lgm} but their performance depend on the existence of algorithms to compute accurately and cheaply the exponential and related functions of a matrix, or their action on vectors.

For example, the solution of the linear equation
\[
   {\bf u}' = A {\bf u} + {\bf b}, \qquad {\bf u}(0)={\bf u}_0
\]
${\bf u},{\bf b}\in \C^{d}$,
$A\in \C^{d\times d}$,
can be written as
\[
  {\bf u}(t) = \e^{tA}{\bf u}_0 + t\varphi_1(tA) {\bf b}
\]
where $\varphi(x)=\frac{\e^{x}-1}{x}$ which requires to compute either the functions of matrices $\e^{tA}$ and $\varphi_1(tA)$ or their action on vectors, being the best choice depending on the particular problem. The solution can also be written in different forms involving only the exponential or the $\varphi_1$ function like, for example
\[
  {\bf u}(t) = \e^{tA}{\bf u}_0 + A^{-1}(\e^{tA} {\bf b} - {\bf b}) 
     = {\bf u}_0 + t\varphi_1(tA) (A{\bf u}_0 + {\bf b})
\]
the first one requiring the computation of the inverse of the matrix $A$ and the second one requires an accurate evaluation of $\varphi_1(tA)$ or its action on a vector. Even in this last case, if the scaling and squaring technique is applied and one considers that
\[
 2t \varphi_1(2tA) = (\e^{tA}+I) \varphi_1(tA)
\]
then, both functions $\e^{tA}$ and $\varphi_1(tA)$ have to be simultaneously evaluated.

In discretised hyperbolic equations frequently one has to solve the equation
\[
  {\bf u}'' + A^2 {\bf u} = 0, \qquad {\bf u}(0)={\bf u}_0, \ {\bf u}'(0)={\bf u}'_0,
\]
with solution
\[
  {\bf u}(t) = \cos(tA) {\bf u}_0 + A^{-1}\sin(tA) {\bf u}'_0
\]
which requires to compute the trigonometric functions of matrices or their actions on vectors. As in the previous case, the scaling and squaring technique can be used with the double angle formulae
\[
  \cos(2A)=2\cos^2(A)-I, \qquad
  \sin(2A)=2\sin(A)\cos(A),
\]
which needs to compute both functions simultaneously.

In some other cases, given a the transition matrix which gives the flow for the evolution of a differential equation, one can be interested to obtain the generator matrix for this problem, and it can be computed with the logarithm of the transition matrix.
Obviously, the particular method to be used will depend on the size and the structure of the matrices whose functions are desired to be computed or their action on vectors.


Then, given a matrix  $A\in \C^{d\times d}$, the goal is to compute $f(A)$ where $f(x)$, with  $x\in \C$, is an analytic function near the origin, e.g. $e^x, \cos(x),\sin(x),\log(1+x)$ or the $\varphi$-functions, i.e. $\varphi_k(x)=(e^x-p_k(x))/x^k$ where $p_k$ is the $k$th order Taylor polynomial to the exponential.

To compute $f(A)$, in this work we consider the case in which 
\[
(I+\alpha A)^{-1} \qquad  \mbox{or}  \qquad (I+\alpha A)^{-1}{\bf v}, 
\]
with $\alpha$ a sufficiently small scalar and ${\bf v}$ a vector, can be efficiently computed. Notice that:

\begin{itemize}
	\item If $A$ is a dense matrix: It is well known that $(I+\alpha A)^{-1}$ can be computed at 4/3 times the cost of the product of two dense matrices. The inverse of a dense matrix can be computed, for example, using an $LU$ decomposition and solving $d$ upper and lower triangular systems with a total number of $\frac83 d^3$ flops or, equivalently, a total cost similar to 4/3 matrix-matrix products.
	\item If $A$ is a large and sparse matrix: the solution of $(I+\alpha A)^{-1}{\bf v}$ can be efficiently carried out in many cases using, for example, incomplete $LU$ or Choleski factorization, the conjugate (or bi-conjugate) gradient method with preconditioners, etc. \cite{golub93mc,trefethen97nla}. For example, if $A$ is tridiagonal (or pentadiagonal) then $(I+\alpha A)$ is also tridiagonal (or pentadiagonal) and, as we will see in more detail, the system $(I+\alpha A){\bf x}={\bf v}$ can be solved with only $8d$ flops ($15d$ flops for pentadiagonal matrices) which can be considered very cheap since the product $A{\bf v}$ already needs $5n$ flops ($9d$ flops for pentadiagonal matrices).

	\item Quantum computation emerged about two decades ago and recently this is a field of enormous research interest. It is claimed that a high performance can be achieved for solving linear algebra problems of large dimension \cite{biamonte17qaf,harrow09qaf}. For instance, in \cite{biamonte17qaf} it si mentioned: {\it The key ingredient behind these methods is that the quantum state
of $n$ quantum bits or qubits is a vector in a $2n$-dimensional complex vector space; performing a quantum logic operations or a measurement
on qubits multiplies the corresponding state vector by $2n\times 2n$ matrices. By building up such matrix transformations, quantum computers
have been shown to perform common linear algebraic operations such as 
Fourier transforms, finding eigenvectors and eigenvalues,
and solving linear sets of equations over $2n$-dimensional vector spaces in time that is polynomial in $n$, exponentially faster than their best known classical counterparts.}
	
\end{itemize}


Two steps are frequently considered when computing most functions, $f(A)$, or their action on vectors:

\begin{itemize}
	\item If the norm of the matrix $A$ is not sufficiently small, an scaling is usually applied that depends on the function to be computed. For example, to compute $\e^A$ one can consider $e^{A/N}$ with $N$ such that the norm of $A/N$ is smaller than a given value and then $e^{A/N}$ is accurately approximated. If $N=2^s$, then $s$ squaring are finally applied. For trigonometric functions, alternative recurrences like the double angle formula can be applied, etc. If, given ${\bf v}\in \C^{d}$, one is interested to compute $e^{A/N}{\bf v}$ the recurrence  ${\bf v}_{n}=e^{A/N}{\bf v}_{n-1}, n=1,2,\ldots,N$, with ${\bf v}_0={\bf v}$ is applied.
 \item One has to compute the scaled function, say $f(B)$ with $B$ depending on $A$, e.g. $B=A/N$, that can be written as a power series expansion
\begin{equation}\label{eq:fB}
   f(B)=\sum_{k=0}^{\infty} a_kB^k.
\end{equation}
 Then, high order rational Chebyshev or Pad\'e approximants, or polynomial approximations are frequently considered to approximate the formal solution following some tricks that allow to carry their computations with a reduced number of operations \cite{almohy09ans,bader19ctm,gallopoulos92eso,sastre19btc}. For example, a Taylor polynomial of degree 18 can be computed with only 5 matrix-matrix products \cite{bader19ctm} or a diagonal Pad\'e approximation, that approximates $e^x$ up to order $x^{26}$, can be computed with only 6 matrix-matrix products and one inverse \cite{almohy09ans}. On the other hand, the computation $f(B){\bf v}$ is frequently carried out using Taylor or Krylov methods \cite{almohy11ctaf,gallopoulos92eso,hochbruck97oks,hochbruck10ei} because the scaling-squaring technique can not be used in this case.
\end{itemize}

In some cases the numerical methods to solve these problems have some stages which can be computed in parallel and this is considered as an extra bonus of the method. However, we are interested on numerical schemes that are built from the very beginning to be used in parallel. The goal of this work is to present a procedure that allows to approximate any function as a linear combination of simple functions that can be evaluated independently so, they can be computed in parallel. In addition, they can be used to approximate simultaneously several functions of matrices too. 
We will also show how similar schemes can be used to compute the action of these functions on vectors.


\subsection{Fractional decomposition}


A technique that has already been used in the literature is the 
approximation of functions by rational approximations
\cite{calvetti95ipf,chehab16pmf,gallopoulos16pim,gallopoulos89otp,gallopoulos92eso,gander13pap,sidje98eas} in which we can apply a fractional decomposition.
Given a function $f(x)$, it can be approximated by a rational function, $r_{n,m}(x)$, such that $r_{n,m}(x)\simeq f(x)$ for a range of values of $x$, where $r_{n,m}(x)=\frac{p_n(x)}{q_m(x)}$, and $p_n(x),q_m(x)$ are polynomials of degree $n$ and $m$, respectively. One can then consider the fractional decomposition
\[
  f(x)\simeq r_{n,m}(x)=\frac{p_n(x)}{q_m(x)} = \sum_{i=1}^m w_i\frac{1}{x-c_i} 
	  + s_{n-m}(x),
\]
where $s_{n-m}(x)$ is a polynomial of degree $n-m$ if $n\geq m$ or $0$ otherwise, and the right hand side can be computed in parallel. If $n-m>2$ the cost can be dominated by the cost to evaluate the polynomial $s_{n-m}(x)$.

The choice of the polynomials $p_n(x), q_m(x)$ depends on the particular methods used, i.e. rational Pad\'e or Chebyshev approximations, and the main trouble is that for most functions of practical interest the roots of $q_m(x)$, i.e. the coefficients $c_i$ are complex making the computational cost about four times more expensive. This can be partially solved if one considers an incomplete fractional decomposition \cite{calvetti95ipf}. Since the complex roots of $q_m(x)$ occur in pairs, one can decompose it as follows 
\[
  f(x)\simeq r_{n,m}(x)=\frac{p_n(x)}{q_m(x)} = 
	\sum_{i=1}^{m_1} w_i\frac{1}{x-c_i} +  
	\sum_{i=m_1+1}^{m_1+m_2} \frac{2Re(w_i)x-2Re(w_ic_i^*)}{x^2-2Re(c_i)x+|c_i|^2} 
	  + s_{n-m}(x),
\]
where $m=m_1+2m_2$ and $c_1,\ldots,c_{m_1}$ are the real roots. Then, some processors have to compute one product, $x^2$, and one inverse which altogether  is nearly twice the cost of one inverse.

Notice that for each function one has to find the polynomials $p_n(x), q_m(x)$ for different values of $n$ and $m$, and then to evaluate the fractional decomposition, making this procedure less attractive

We simplify this procedure using only real coefficients and making the search for the coefficients of the fractional decomposition trivial for most functions as well as to show how to adapt the procedure when the matrix $A$ has different properties.
The paper is organised as follows: Section~\ref{sec:MainIdea} presents the main idea to build the methods. An error analysis is presented in Section~\ref{sec:error}. Section~\ref{sec:examples} illustrates how to build some particular methods which are numerically tested in Section~\ref{sec:NumExamples}. Finally, Section~\ref{sec:conclusions} collects the conclusions as well as future work.

\section{Approximating functions by simple fractions} \label{sec:MainIdea}

The main idea of this work is quite simple, notice that
\[
 F_i(B)= (I-c_iB)^{-1}=\sum_{k=0}^{\infty} c_i^k B^k
\]
for sufficiently small values of $\|c_iB\|$ and whose computational cost, as previously mentioned, can be considered, for general dense matrices, as 4/3 matrix-matrix products. Notice that each function $F_i$, for different values of $c_i$, can be computed in parallel and then, if $P$ processors are available, we can compute
\begin{equation} \label{eq:rat1}
  r_s^{(P)}(x)=\sum_{i=1}^P b_iF_i(x) 
\end{equation}
where
\begin{equation} \label{eq:approx1}
  f(x)-r_s^{(P)}(x)= {\cal O}(x^{s+1})
\end{equation}
with $s=P-1$
if the coefficients are chosen such that $c_i\neq c_j$ for $i\neq j$ and the coefficients $b_i$ solve the following simple linear system of equations
\begin{equation} \label{eq:LinearEqsbi}
   \sum_{i=1}^P b_ic_i^k=a_k, \qquad k=0,1,\ldots,P-1,
\end{equation}
where the coefficients $a_k$ are known from \eqref{eq:fB}.
Obviously, once the functions $F_i$ are computed, different functions $f(B)$  (with different values for  the coefficients $a_k$) can be simultaneously approximated just by
looking for a new set of coefficients $b_i$ in \eqref{eq:rat1}, say  $\hat b_i$, such the corresponding equations \eqref{eq:LinearEqsbi} with the new coefficients $a_k$ are satisfied.
It remains the problem on how to choose the best set of coefficients, $c_i$, for each function, and this will depend on the number of processors available, $P$, the accuracy or stability desired, etc. If we choose $s+1<P$ then $P-s-1$ coefficients $b_i$ can be taken as free parameters for optimization purposes.

 For example, the 4th-order diagonal Pad\'e approximation to the exponential is given by 
\begin{eqnarray} \label{eq.Pade4}
 r_{2,2}(x)&=&\frac{1+\frac12 x + \frac1{12}x^2}{1-\frac12 x + \frac1{12}x^2}=  
  1+\frac{6+i\, 6\sqrt{6}}{x-(3-i\, \sqrt{3})}+\frac{6-i\, 6\sqrt{6}}{x-(3+i\, \sqrt{3})} \\
	&=& 
	1+2 \mbox{Re}\left( \frac{6+i\, 6\sqrt{6}}{x-(3-i\, \sqrt{3})} \right) \nonumber 
\end{eqnarray}
with $r_{2,2}(x)=e^x+{\cal O}(x^5)$, but complex arithmetic is involved, and  $r_{2,2}(x)$ is only valid to approximate the exponential function. Obviously, $r_{2,2}(x)$ can also be computed at the cost of one product, $x^2$, plus one inverse, $(1-x/2 + x^2/12)^{-1}$. On the other hand, a 4th-order Pad\'e  approximation to the function $\varphi_1(x)$ is given by
\begin{equation} \label{eq.Pade4phi1}
 \tilde r_{2,2}(x)=\frac{1+\frac1{10} x + \frac1{60}x^2}{1-\frac25 x + \frac1{20}x^2}= 
\frac13+2 \mbox{Re}\left( \frac{\frac73+i\, 8}{x-(4-i\, 2)} \right)
\end{equation}
with $\tilde r_{2,2}(x)=\varphi_1(x)+{\cal O}(x^5)$ which has different complex roots.

However, if five processors are available we can take, for example
\[
  r_4^{(5)}(x)=\sum_{i=1}^5 b_i\frac{1}{1-c_ix},
\]
where the coefficients $c_i$ can be chosen to optimise the performance of the method. One set of coefficients $c_i$ can be optimal to get accurate results
while other set of coefficients can be more appropriate when, for example, the matrix $B$ is positive or negative definite 
(the coefficients $c_i$ with appropriate size and sign can be chosen to optimize stability). 
If we take $c_1=0$ and $b_1=a_0-(b_2+b_3+b_4+b_5)$ then only four processors are required.
 Once the values for $c_i$ are fixed, the coefficients $b_i$ are trivially obtained. 

For example, if we choose $c_i=1/(i+1), \ i=1,\ldots,5$ then the system \eqref{eq:LinearEqsbi} with $a_k=1/k!, \ k=0,1,\ldots,4$ has the solution
\begin{equation} \label{eq.Frac4bi0}
  b_1=\frac{1}{3}, \quad
  b_2=-18, \quad
  b_3=128, \quad
  b_4=-\frac{625}{3}, \quad
  b_5=99,
\end{equation}
which corresponds to an approximation to the exponential that is already more accurate than the previous 4th-order Pad\'e approximation \eqref{eq.Pade4}.
In addition, if one takes $a_k=1/(k+1)!$, then we will obtain the solution
\[
  b_1=\frac{7}{18}, \quad
  b_2=-9, \quad
  b_3=\frac{128}{3}, \quad
  b_4=-\frac{500}{9}, \quad
  b_5=\frac{45}{2}, \quad
\]
which corresponds to a 4th-order approximation to the function $\varphi_1(x)$, 
or a 4th-order approximation to the function $\log(1-x)$ for sufficiently small $x$ can be obtained if we take
\[
  b_1=-\frac{35}{3}, \quad
  b_2=\frac{153}{2}, \quad
  b_3=-160, \quad
  b_4=\frac{625}{6}, \quad
  b_5=-9, \quad
\]
although this set of coefficients $c_i$ is not necessarily the optimal one for these functions.

If one is interested to compute the action of the function on a vector one has to consider
\[
  f(B){\bf v} \simeq \sum_{i=1}^P b_i {\bf v}_i, \qquad \mbox{where} \qquad 
	(I-c_iB){\bf v}_i={\bf v}
\]
which, as previously, can be computed by solving linear systems of equations in each processor in parallel.

\subsection{Functions of tridiagonal matrices acting on vectors}
Just as an illustration of functions of sparse matrices acting on vectors, let us consider the computation $f(B){\bf v}$ where $B$ is a tridiagonal matrix. The product  $B{\bf v}$ is done with $5d$ flops so, a Krylov or Taylor polynomial of degree $K$, requires $5d\times K$ flops, in addition to the evaluation of the function $f$ for a matrix of dimension $K\times K$ for the Krylov methods. However, since $(I-c_iB)$ is also tridiagonal then  $(I-c_iB){\bf v}_{i}={\bf v}$ can be solved with only $8d$ flops (using, for example, the Thomas algorithm). Then, for problems of relatively large size this is a significant saving in the computational cost (if the cost to communicate between processors can be neglected or is not dominant). 
 This problem will be considered later with more detail.
Linear systems for tridiagonal matrices can also be solved in parallel (see, e.g. \cite{press97nr}), but this would require a second level of parallelism which is not considered in this work.

\section{Error analysis}\label{sec:error}

Forward and backward error analysis can be easily carried out for the new classes of approximations. For example, if one is interested in forward error bounds, we have that
\[
  \|f(B)-r_s^{(P)}(B)\| =\|\sum_{k=0}^{\infty} (a_k-\alpha_k)B^k\|\leq
	\sum_{k=0}^{\infty} |a_k-\alpha_k|\|B\|^k
\] 
where $\alpha_k=\sum_{i=1}^Pb_ic_i^k$. Since, in general, the coefficients $a_k,b_i,c_i$ are known with any desired accuracy, one can take a sufficiently large number of terms in the summation to compute the error bounds with several significant digits. For different values of the tolerance, $\epsilon$, we can find values, say $\theta_s$, such that if $\|B\|<\theta_s$ the error is below  $\epsilon$. Obviously, sharper error bounds can be obtained if additional information of the matrix is known, e.g. when the bounds are written in terms of  $\|B^{k}\|^{1/k}$ for some values of $k$ (see, e.g. \cite{almohy09ans}).

Alternatively, error bounds from the backward error analysis can be obtained if one considers that 
$r_s^{(P)}(B)$ can formally be written as
\[
  r_s^{(P)}(B) = f(B+\Delta B)
\] 
where 
\[
  \Delta B= f^{-1}(r_s^{(P)}(B))-B = h_s(B) = \sum_{k=0}^{\infty} \beta_k B^k
\]
where $\beta_k=0$ for $k=0,1,\ldots,s$ since $r_s^{(P)}(x)$ coincides with the Taylor expansion of $f(B)$ up to order $s$, but this is not necessarily the case in schemes like Chebyshev approximations where $s=0$. Then, it is clear that $\|h_P(B)\|\leq \tilde h_s(\|B\|)$ where
\[
   \tilde h_P(x)=\sum_{k=0}^{\infty} |\beta_k| x^k, 
\]
and thus
\[
  \frac{\|\Delta B\|}{\|B\|} =
  \frac{\|h_P(B)\|}{\|B\|} \leq
  \frac{\tilde h_P(\|B\|)}{\|B\|}.
\]
These previous error analysis can be easily carried for different functions $f(B)$ and for different choices of the coefficients $c_i$ and their associated values for the coefficients $b_i$.

However, one has to take into account that the coefficients $b_i$ will strongly depend on the choice of the coefficients $c_i$. On one side, we have that taking very small values for the coefficients $c_i$ usually allows to apply the methods for matrices with large norms, but then the solution of the linear system of equations \eqref{eq:LinearEqsbi} will have, in general, large values for the coefficients $b_i$, and this can cause badly conditioned methods when high accuracy (say, near to round off errors) is desired. This can be partially solved taking additional processors or slightly modifying the methods as we will see.

\section{Illustrative examples: Exponential and phi-functions}\label{sec:examples}

We illustrate how to build and optimise some approximations of the matrix exponential and the $\varphi_1$-function taking into account the error analysis when only a small number of processors are available and relatively low order methods are considered.
For simplicity in the presentation we will only consider forward error bounds.

\subsection{4th-order approximations}

One of the best 4th-order methods to approximate the the matrix exponential is the $r_{2,2}(x)$ Pad\'e approximation which in a serial algorithm needs one product and one inverse and similarly for $\tilde r_{2,2}(x)$ to approximate the $\varphi_1$-function. A forward error bound is given by
\[
  \|e^{B}-r_{2,2}(B)\|\leq \epsilon_{4,pad}(\|B\|), \qquad 
  \|\varphi_1(B)-\tilde r_{2,2}(B)\|\leq \tilde \epsilon_{4,pad}(\|B\|)
\]
where
\[
\epsilon_{4,pad}(x) = \sum_{k=5}^{\infty} \left| \frac1{k!}-d_k   \right| x^k, \qquad \quad
\tilde \epsilon_{4,pad}(x) = \sum_{k=5}^{\infty}  \left| \frac1{(k+1)!}-\tilde d_k   \right| x^k.
\]
Here, $d_k$ and $\tilde d_k$ are the coefficients from the Taylor expansion of $r_{2,2}(x)$ and $\tilde r_{2,2}(x)$, respectively.

Let us now consider, for simplicity in the search of an optimized parallel method, the following choice for the coefficients $c_i$ for an scheme using four processors: 
\begin{equation} \label{eq.Frac4ci}
  c_1=0, \quad 
	c_{2}=\frac1{\alpha}, \quad 
	c_{3}=-\frac1{\alpha}, \quad 
	c_{4}=\frac1{2\alpha},  \quad 
	c_{5}=-\frac1{2\alpha},
\end{equation}
where the coefficients $b_i$ are trivially obtained. Then, the proposed method to be used with four processors is written in terms of one free parameter, $\alpha$. 
The forward error bound is given by
\[
  \|e^{B}-r_{4}^{(4)}(B)\|\leq \epsilon_{4}(\|B\|), \qquad
  \|\varphi_1(B)-\tilde r_{4}^{(4)}(B)\|\leq \tilde \epsilon_{4}(\|B\|)
\]
where
\[
\epsilon_{4}(x) = \sum_{k=5}^{\infty} \left| \frac1{k!}-d_k(\alpha)   \right| x^k, \qquad \qquad
\tilde \epsilon_{4}(x) = \sum_{k=5}^{\infty} \left| \frac1{(k+1)!}-\tilde d_k(\alpha)   \right| x^k.
\]
Here, $d_k(\alpha)=\sum_{i=2}^5b_ic_i^k$ and $\tilde d_k(\alpha)=\sum_{i=2}^5\tilde b_ic_i^k$ are the coefficients from the Taylor expansion of $r_{4}^{(4)}(x)$ and $\tilde r_{4}^{(4)}(x)$, respectively, which depend on $\alpha$. A simple search shows that the optimal solution which minimises $\epsilon_{4}(x)$ for most values of $x$ about the origin occurs (approximately) for $\alpha=5$, and then
\begin{equation} \label{eq.Frac4bi}
  b_1= \frac{128}{3},\quad
  b_2= \frac{85}{3},\quad
  b_3= \frac{20}{9},\quad
  b_4= -\frac{515}{9},\quad
  b_5= -15.
\end{equation}
For this choice of $\alpha$ we can still get an approximation for the $\varphi_1$ function that is nearly as accurate as the previous Pad\'e approximation. The optimal choice for this function is however (approximately) for $\alpha=6$, and then
\[
  b_1= \frac{71}{5},\quad
  b_2= \frac{117}{10},\quad
  b_3= \frac{7}{10},\quad
  b_4= -\frac{104}{5},\quad
  b_5= -\frac{24}{5},
\]
which gives more accurate results than the Pad\'e scheme.

\begin{figure}[h!]
  \begin{center}
       \includegraphics[scale=0.5]{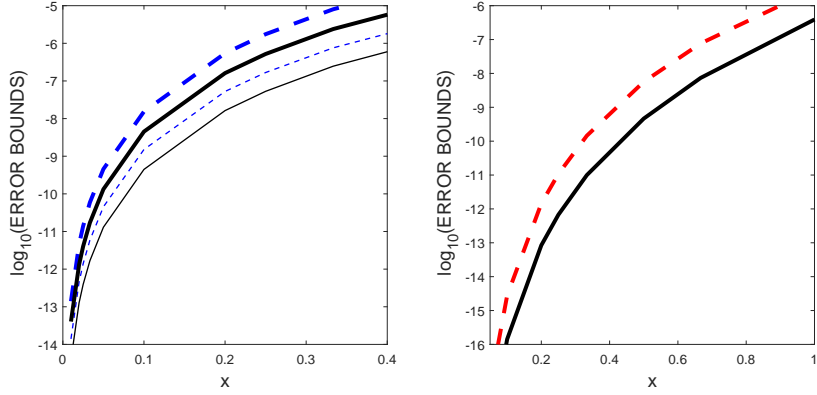}
      \caption{Left: Error bounds for the 4th-order approximations to the exponential function (thick lines) and to the $\varphi_1$-function (thin lines) using Pad\'e approximations (dashed lines) and the fractional approximations (solid lines). Right:  Error bounds for the 8th-order approximations to the exponential function using the Taylor approximation (dashed line) and the fractional approximation (solid line).}
    \label{Fig.0}
  \end{center}
\end{figure}

In Fig.~\ref{Fig.0}, left panel, we show the values of the different error bounds, $\epsilon_{4,pad}(x)$, $\tilde \epsilon_{4,pad}(x)$, $\epsilon_{4}(x)$ and $\tilde \epsilon_{4}(x)$, versus $x$. The thick lines correspond to the 4th-order approximations to the exponential function  and the thin lines correspond to the $\varphi_1$-function when using Pad\'e approximations (dashed lines) and the fractional approximations (solid lines).
 Then, given a tolerance, $tol$, one can easily find the value $\theta$ such that $\epsilon(x)<tol$ for $x<\theta$, i.e. $\epsilon(\|B\|)<tol$ for $\|B\|<\theta$ where $\epsilon$ denotes the desired error bound. It is clear that the new methods are more accurate for all values of $\|B\|$ and are also faster to be computed when done in parallel, and we have to remark that the new methods are not fully optimized.

\subsection{8th-order approximations}

In a serial computer, the 8th-order Taylor approximation can be computed with only 3 matrix-matrix products as follows \cite{bader19ctm}
\begin{equation}  \label{Algorithm83eB}
\begin{aligned}
   B_2 & = B^2,\\
   B_4 & = B_2(x_1 B+x_2 B_2),\\
   B_8 & = (x_3 B_2+B_4)(x_4I+x_5B+x_6 B_2+x_7B_4),  \\
   T_8(B)  & = y_0 I + y_1 B +y_{2} B_2+ B_8,
\end{aligned}  
\end{equation}
with
\[
\begin{array}{llll}
x_1=\displaystyle x_3\frac{ 1 + \sqrt{177}}{88},& 
x_2= \displaystyle \frac{1 + \sqrt{177}}{352}{x_3},&
x_4= \displaystyle \frac{-271 + 29\sqrt{177}}{315 x_3}, \\
x_5= \displaystyle \frac{11 (-1 + \sqrt{177})}{1260 x_3},&
x_6=  \displaystyle \frac{11 (-9 + \sqrt{177})}{5040 x_3},&
x_7=  \displaystyle \frac{89 - \sqrt{177}}{5040 x_3^2},&\\
y_0=1,& 
y_1= 1,& 
y_2 = \displaystyle \frac{857 - 58\sqrt{177}}{630}, \\
x_3 = 2/3,
\end{array}
\]
being the most efficient scheme for this order of accuracy.
In spite the coefficients $x_i,y_i$ are not rational numbers one can check that $T_8$ is exactly the Taylor expansion to order 8, and a similar method can be obtained for the $\varphi_1$-function which we do not consider for simplicity.
The forward error bound is given by
\[
  \|e^{B}-T_{8}(B)\|\leq \epsilon_{8}^T(\|B\|),
	\qquad  \mbox{where} \qquad
\epsilon_{8}^T(x) = \sum_{k=9}^{\infty} \frac1{k!}x^k.
\]

Let us take, for example 
\[
  r_8^{(8)}=b_1+\sum_{i=2}^{9}b_i\frac1{1-c_ix}
\]
(with $c_1=0$) and 
\begin{eqnarray}
 & & 
  c_{2}=\frac1{\alpha}, \quad 
	c_{3}=-\frac1{\alpha}, \quad 
	c_{4}=\frac2{3\alpha},  \quad 
	c_{5}=-\frac2{3\alpha},  \quad \nonumber \\
	& & c_{6}=\frac1{2\alpha}, \quad 
	c_{7}=-\frac1{2\alpha}, \quad 
	c_{8}=\frac2{5\alpha},  \quad 
	c_{9}=-\frac2{5\alpha},  \label{eq.Frac8ci}
\end{eqnarray}
to be used with eight processors and written in terms of one free parameter, $\alpha$. The optimal value for this choice of coefficients $c_i$ which minimize the forward error bound, $\epsilon_8(x)$ corresponds (approximately) to $\alpha=5$, and then 

\begin{equation}\label{eq.Frac8bi}
	\begin{array}{lll}
\displaystyle b_1= -\frac{9979069}{32256},\qquad 
& \displaystyle b_2= -\frac{1995521}{254016},\quad 
& \displaystyle b_3= -\frac{392009}{254016}, \nonumber\\
 & & \\
\displaystyle b_4= \frac{520866369}{802816},\quad 
& \displaystyle b_5= \frac{48898161}{802816}, 
& \displaystyle b_6= -\frac{26686735}{11907},\nonumber\\
 & & \\
\displaystyle b_7= -\frac{3892615}{11907},\quad 
& \displaystyle b_8= \frac{353067578125}{195084288},\quad 
& \displaystyle b_9= \frac{71873828125}{195084288}. 
\end{array}
\end{equation}

In the right panel of Fig.~\ref{Fig.0}  we show the values of $ \epsilon_{8}^T(x)$ (dashed line) and $\epsilon_{8}(x)$ for $\alpha=5$ (solid line) for different values of $x$. We observe that this method, with such a simple optimization search, already provides more accurate results than the Taylor approximation being up to slightly more than twice faster.

\subsection{High order approximations}

From the error bounds analysis we can deduce that when high accuracy is desired, it is usually more efficient to consider higher order approximations rather than taking a larger value of $N$ in the approximation $B=A/N$ and next to consider the recurrence like the squaring. This usually occurs up to a relatively high order. 

However, in the construction of different methods we have noticed that in the optimization process for the fractional approximation, better error bounds are obtained when small values of the coefficients $c_i$ are taken and this usually leads to large values for the coefficients $b_i$ whose absolute values typically grow with the order of accuracy (the coefficients $b_i$ of our 8th-order approximation are considerably larger in absolute value than the corresponding coefficients for the 4th-order method). 

Then, if nearly round off accuracy is desired, e.g. double precision, the new schemes could not be well conditioned and can suffer from large round off errors. This problem can be partially reduced in different ways. 
If the computation of $(I-c_iB)^{-1}I$ has similar cost as  $(I-c_iB)^{-1}B$ we can take into account that $\frac1{1-c_ix}=1+ \frac{c_ix}{1-c_ix}$ and we have that
\begin{equation} \label{eq.ParallRoundOff}
  r_s^{(P)}(x)=a_0 + \sum_{i=1}^P b_i\frac{c_ix}{1-c_ix}
\end{equation}
where we have considered that $\sum_{i=1}^Pb_i=a_0$. This decomposition usually has smaller round-off errors for relatively small values of $|c_ix|$ which is usually the case when close to round off errors are desired.

Alternatively, if there is no restriction on the number of processors, one can consider in \eqref{eq:approx1} $P>s+1$ allowing for $P-s-1$ free coefficients $b_i$ that can be chosen jointly with appropriate values for the $c_i$'s to reduce round off errors. We can end up with an optimisation problem with many free parameters (the coefficients $c_i$ in addition to $P-s-1$ coefficients $b_i$) and we leave this as an interesting open problem to be analysed for different functions of matrices, number of processors used and choices of the order of accuracy. In \cite{jarlebring21cgf} the authors show how to optimise the search of the coefficients in some approximations of functions of matrices that hopefully could be used for this problem too.

There are also hybrid methods that could be used to reduce round off errors. For example, we can consider
\begin{equation} \label{eq.ParallHybrid}
  r_s^{(P)}(x)=d_0+d_1x+d_2x^2+\sum_{i=1}^{P-1} b_i\frac{1}{1-c_ix},
\end{equation}
which can be considered as a generalisation of the fractional decomposition of $r_{n+2,n}$ rational Pad\'e approximations. Here, $d_0+d_1x+d_2x^2$ can be computed in one processor (at the cost of one matrix-matrix product, that is slightly cheaper than the inverse of a matrix) and the partial fractions are computed in the remaining $P-1$ processors. The coefficients $b_i$ must satisfy (for $P\geq 2$)
\[
   \sum_{i=1}^{P-1} b_ic_i^k=a_k, \qquad k=3,4,\ldots,P+1
\]
and 
\[
   d_k=a_k - \sum_{i=1}^{P-1} b_ic_i^k, \qquad k=0,1,2,
\]
which allows to get methods up to order $s=P+1$ with $P$ processors and, hopefully, smaller coefficients $b_i$ in absolute value. 
%
It is then expected that a more careful search will lead to new schemes with considerably reduced round off errors.

We illustrate this procedure in the search of some 10th-order methods for the matrix exponential.

\subsubsection{10th-order approximations to the exponential}

The 10th-order diagonal Pad\'e approximation which is given by \cite{higham05tsa} and is used as one of the methods implemented in the function \texttt{expm} of {\sc Matlab} is:
\begin{equation}\label{eq.pade10}
	(-u_{5}+v_{5}) r_{5,5}(A) = (u_{5}+v_{5}), \quad 	
	u_{5}  =  A[b_5A_4 +  b_3A_2 + b_1 I ], \quad
	v_{5} =  b_4A_4 + b_2A_2 + b_0 I,
\end{equation}
where $A_2=A^2, A_4=A_2^2, \ (b_0,b_1,b_2,b_3,b_4,b_5)=
(1,\frac{1}{2} ,\frac{1}{9} ,\frac{1}{72},\frac{1}{1008},\frac{1}{30240})$.
It requires 3 products and one inverse (approximately three times more expensive than the inverse to be evaluated by each processor on a parallel method). Note that the Pad\'e approximation for the $\varphi_1$-function has not the same symmetry for the numerator and denominator and has to be computed with four products and one inverse, i.e. four times more expensive that one inverse.

An approximation by using the proposed fractional methods to order ten requires 11 processors (10 processors if one takes, e.g. $c_1=0$ or 9 processors if one considers \eqref{eq.ParallHybrid}). Suppose we have two extra free coefficients $b_i$ in the composition 
 \eqref{eq.ParallHybrid} (11 processors in total) to have some freedom to improve the accuracy of the method (this is just a very simple illustrative example that has not been exhaustively optimised)
\begin{equation}\label{eq.Frac10}
  r_{10}^{(11)}(x)=d_0+d_1x+d_2x^2+\sum_{i=1}^{10} b_i\frac{1}{1-c_ix},
\end{equation}
where, given some values for the coefficients $c_i$ we take as free parameters $b_9,b_{10}$ and solve the following linear systems of equations for the remaining coefficients:
\[
   \sum_{i=1}^{10} b_ic_i^k=\frac1{k!}, \qquad k=3,4,\ldots,10
\]
which provides $b_1,\ldots,b_8$, and 
\[
   d_k=\frac1{k!} - \sum_{i=1}^{10} b_ic_i^k, \qquad k=0,1,2.
\]

We have taken 
\begin{equation}\label{eq.Frac10ciNO}
c_i=\frac1{6+i}, \ \ i=1,2,\ldots,10;\qquad  
b_{9}= -50000, \quad  b_{10}=350000, 
\end{equation}
being the remaining coefficients
\begin{align}\label{eq.Frac10biNO}
d_1&=-\frac{34244933704346617}{14676708556800},  
&d_2&=-\frac{18541933026870559}{428070666240000}, \\
d_3&=-\frac{6798371473106351}{15410543984640000}, & & \nonumber\\ 
b_1&=-\frac{2385751622325187262153}{1585084524134400000},
&b_2&= \frac{3752603696192081}{82668600000},  \nonumber  \\ 
b_3&=-\frac{87230538798639033213}{187904819200000},
&b_4&= \frac{14789110319838821875}{6934744793088}, \nonumber  \\  
b_5&=-\frac{10558563753676365149296981}{2219118333788160000}, 
&b_6&=\frac{3520134037769971}{716800000},  \nonumber \\
b_7&=-\frac{2263057299714115181019509}{1585084524134400000}, 
&b_8&= -\frac{2275831141003773874927}{3095868211200000}.  \nonumber 
\end{align}
This method has an error bound smaller than the error bound provided by the Pad\'e approximation by more than one order of magnitude while being up to three times cheaper in a parallel computer. However, note that
\[
   \max\{|b_i|,|d_i|\}\simeq 4.9\cdot 10^6
\]
so, one can expect large cancellations leading to large round off errors.

This problem can be partially solved by choosing different values for the coefficients $c_i$. Note that we have much freedom in this choice. 
For example, if we take the following values: 
\begin{equation}\label{eq.Frac10ci}
c_{2i-1}=-\frac1{6+2i},\ c_{2i}=\frac1{6+2i}, \quad  i=1,2,3,4,5, \qquad
b_{9}= 2000, \quad b_{10}=-3500, 
\end{equation}
then we get
\begin{align}\label{eq.Frac10bi}
d_1& =-\frac{781562376863}{94371840},\  
&d_2& =-\frac{54849495983}{330301440},\
&d_3& =-\frac{8034429391}{587202560},  \\ 
b_1& =-\frac{57383239}{760320},\quad 
&b_2& =-\frac{115498838729}{239500800},\  
&b_3& = \frac{1648441938671875}{1255673954304},\nonumber \\
b_4& = \frac{56790060546875}{4227858432},\   
&b_5& =-\frac{1476772203681}{298188800},\ 
&b_6& =-\frac{31012455666807}{656015360},\nonumber   \\
b_7& = \frac{4891212112962371}{1295536619520},\ 
&b_8& = \frac{171190903245297593}{3886609858560},\ &\nonumber
\end{align}
where
\[
   \max\{|b_i|,|d_i|\}\simeq 4.7\cdot 10^4.
\]
The forward error bound is slightly worst for this method but, as we will observe in the numerical experiments, round off errors are reduced by nearly three orders of magnitude. 

Obviously, a deeper search taking into account the freedom in the choice of the coefficients $c_i$ as well as the possibility to take more processors with more free parameters $b_i$ and $c_i$ should allow us to considerably reduce both the error bounds as well as the round off error.



\section{Numerical examples}\label{sec:NumExamples}

Given a matrix, $A$, to compute a function, $f(A)$, most algorithms first evaluate a bound to a norm, say $\|A\|$ or $\|A^k\|^{1/k}$ for some values of $k$ and then, according to its value and a tabulated set of values obtained from the error bound analysis, it is chosen the method to be used that gives a result with (hopefully) an error below the desired tolerance. 

To test the algorithms on particular problems allows us to check if the error bounds are sufficiently sharp as well as to observe how big the undesirable round off errors are.

In the numerical experiments we only compute approximations to $\e^{hA}$ or its action on vectors for different values of $h$ using the following methods
\begin{itemize}
	\item $r_{2,2}$: the 4th-order Pad\'e approximation \eqref{eq.Pade4}. Cost: one product and one inverse.
	\item $T_8$: the 8th-order Taylor approximation, $T_{8}(x)$ \eqref{Algorithm83eB}. Cost: three products.
	\item $T_{5}, T_{10}$: the 5th- and 10th-order Taylor approximations used for the action of the exponential on a vector. Cost: five and ten vector-matrix products, respectively.
	\item $r_{5,5}$: the 10th-order Pad\'e approximation \eqref{eq.pade10}. Cost: three products and one inverse.
	\item $R_4$: the 4th-order rational approximation
	\[
  r_4^{(4)}(x)=\sum_{i=1}^5 b_i\frac{1}{1-c_ix},
\]
with coefficients given in \eqref{eq.Frac4ci} and  \eqref{eq.Frac4bi} for $\alpha=5$, designed to be used with  four processors.
	\item $R_5$: the 5th-order rational approximation, $r_5^{(5)}(x)$,
with coefficients given in \eqref{eq.Frac5bi0}, designed to be used with five processors.
	\item $R_8$: the 8th-order rational approximation
	\[
  r_8^{(8)}(x)=\sum_{i=1}^9 b_i\frac{1}{1-c_ix},
\]
with coefficients given in \eqref{eq.Frac8ci} for $\alpha=5$, designed to be used with eight processors.
	\item $R_{10}^*$: the 10th-order rational approximation \eqref{eq.Frac10} with coefficients given in \eqref{eq.Frac10ciNO} and \eqref{eq.Frac10biNO}, designed to be used with eleven processors.
	\item $R_{10}$: the 10th-order rational approximation \eqref{eq.Frac10} with coefficients given in \eqref{eq.Frac10ci} and \eqref{eq.Frac10bi}, designed to be used with eleven processors.
\end{itemize}

As a first numerical test we take 
\[
   A=\texttt{randn(100)}
\] 
i.e. a $100\times 100$ matrix whose elements are normally distributed random numbers and we compute $\e^{hA}$ for different values of $h$. We have done this numerical experiment repeatedly many times with very similar results.

\begin{figure}[h!]
  \begin{center}
       \includegraphics[scale=0.4]{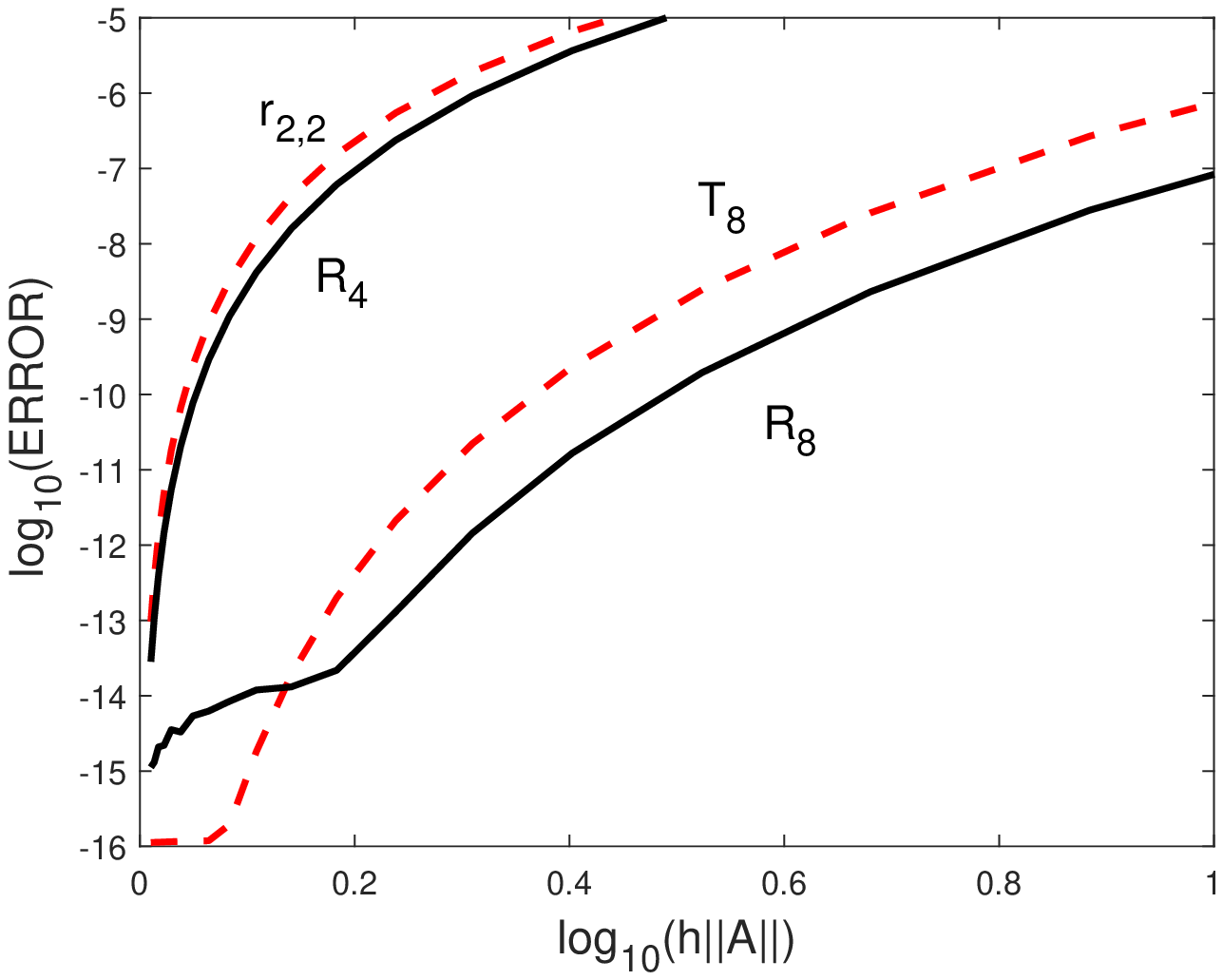}
       \includegraphics[scale=0.4]{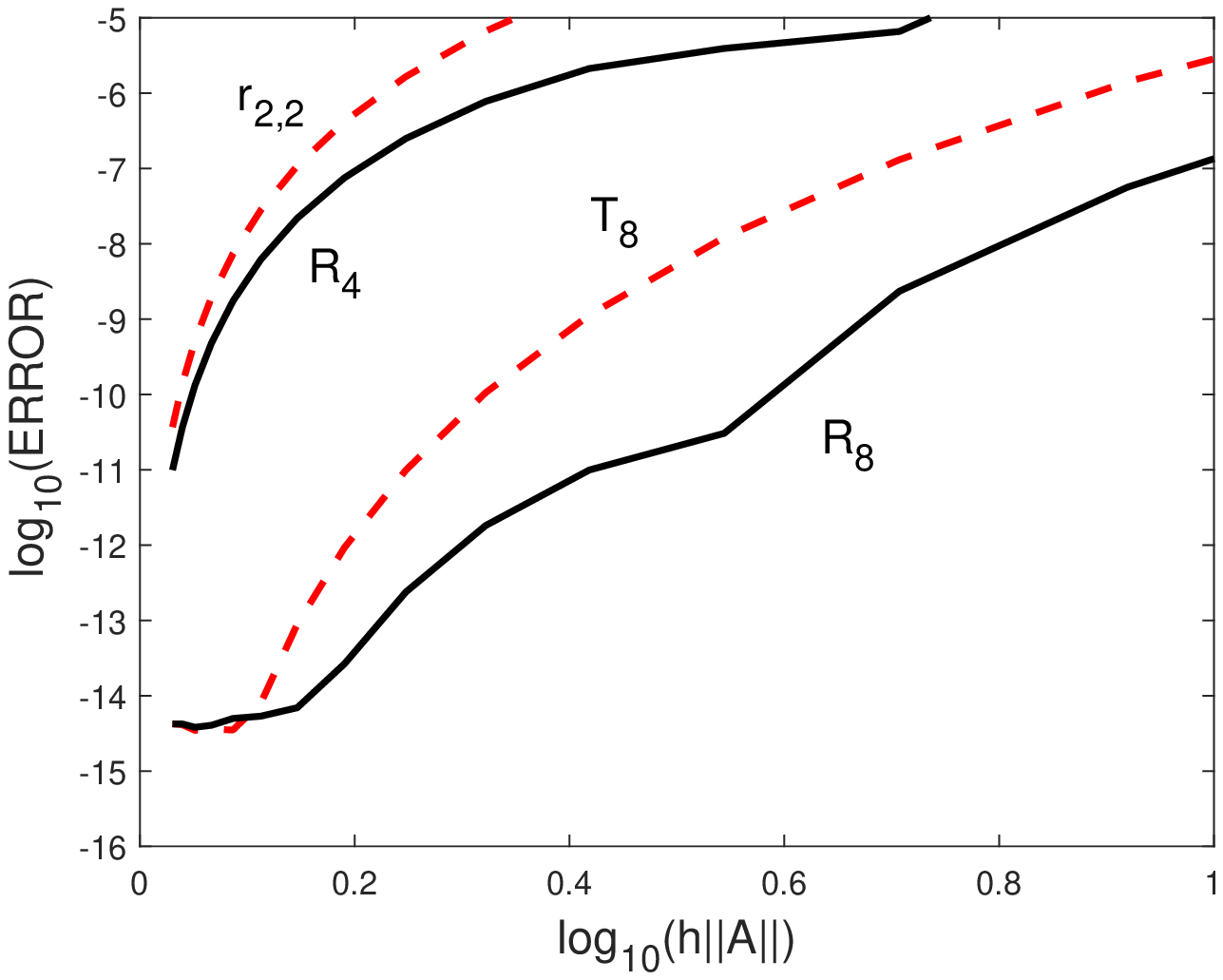}
      \caption{Left:Two-norm error (in logarithmic scale) to compute $\e^{hA}$ with $A$ a random matrix of dimension $d=100$ versus $h\|A\|$. Dashed lines correspond to the results obtained with the 4th-order diagonal Pad\'e method and the 8th-order Taylor method and solid lines are obtained with the new rational methods of the same order. Right: The same with the matrix \eqref{eq.Ex2b} of dimension $d=100$.}
    \label{Fig.2}
  \end{center}
\end{figure}

 Figure~\ref{Fig.2} shows in the left panel the two-norm error (the exact solution is computed numerically with sufficiently high accuracy) versus $h\|A\|$ for these methods. Dashed lines correspond to the results obtained with the diagonal Pad\'e and Taylor methods and solid lines are obtained with the new rational methods. We observe that the new methods have in all cases  similar or higher accuracy than the methods of the same order and built to be computed in a sequential algorithm. There is only a minor drawback in the 8th-order rational method when high accuracy is desired due to round off errors.

Note that the computational cost is not taken into account in the figures since this depends on how efficiently is implemented the rational method in a parallel computer, where the new methods can be up to several times faster to compute than the methods designed for serial algorithms.

We have repeated the numerical experiments with another $100\times 100$ symmetric dense matrix with elements
\begin{equation}\label{eq.Ex2b}
   A_{i,j}=\frac1{1+(i-j)^2}
\end{equation}
and the results are shown in the right panel in  Figure~\ref{Fig.2}. We observe that for this problem round off errors are similar for both 8th-order methods. We have repeated these numerical experiments for different dimensions of the matrices and other matrices with different structures and the results are, in general, qualitatively similar. Notice the results are in agreement with the error bounds shown in Figure~\ref{Fig.0}.

\begin{figure}[h!]
  \begin{center}
       \includegraphics[scale=0.4]{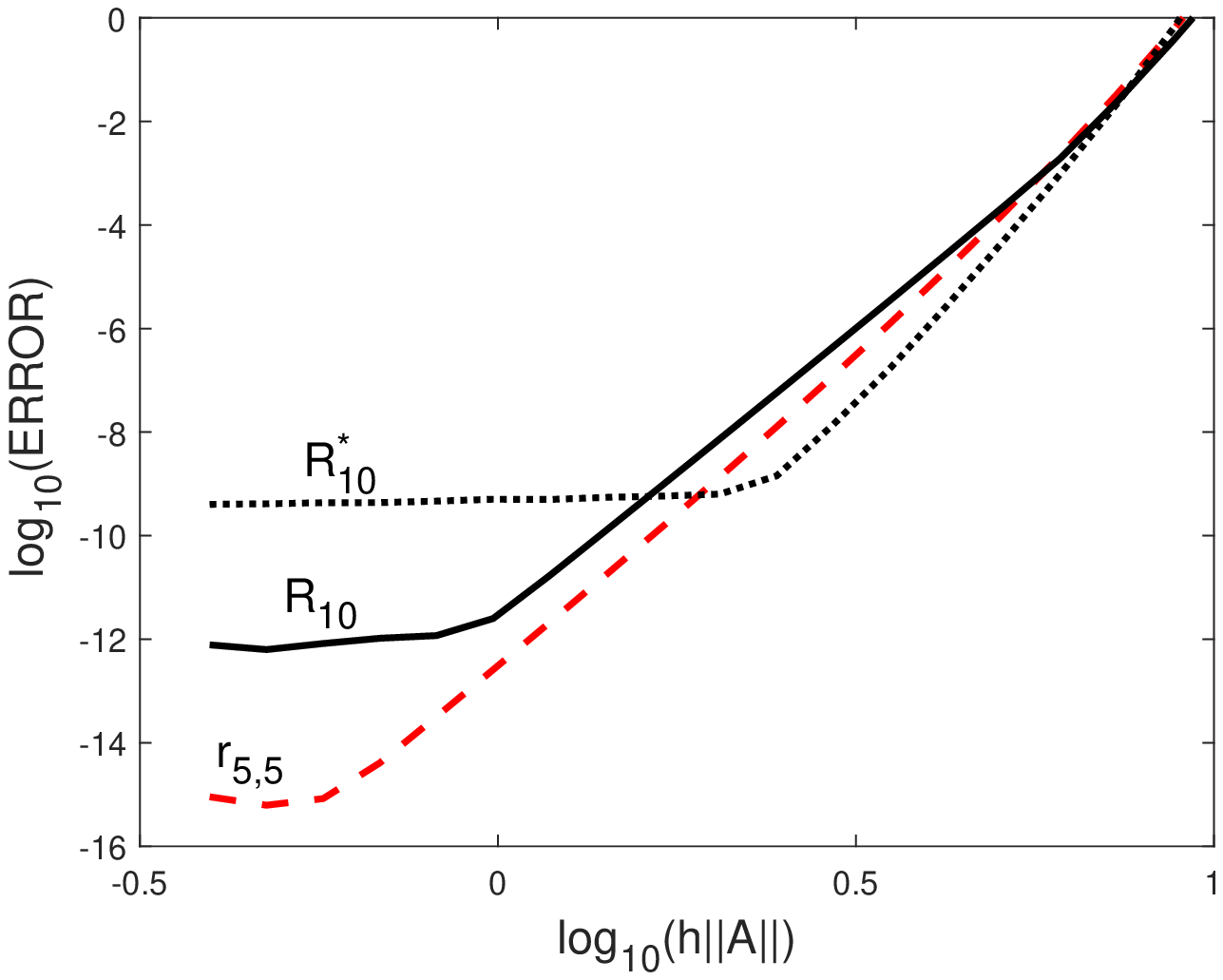}
       \includegraphics[scale=0.4]{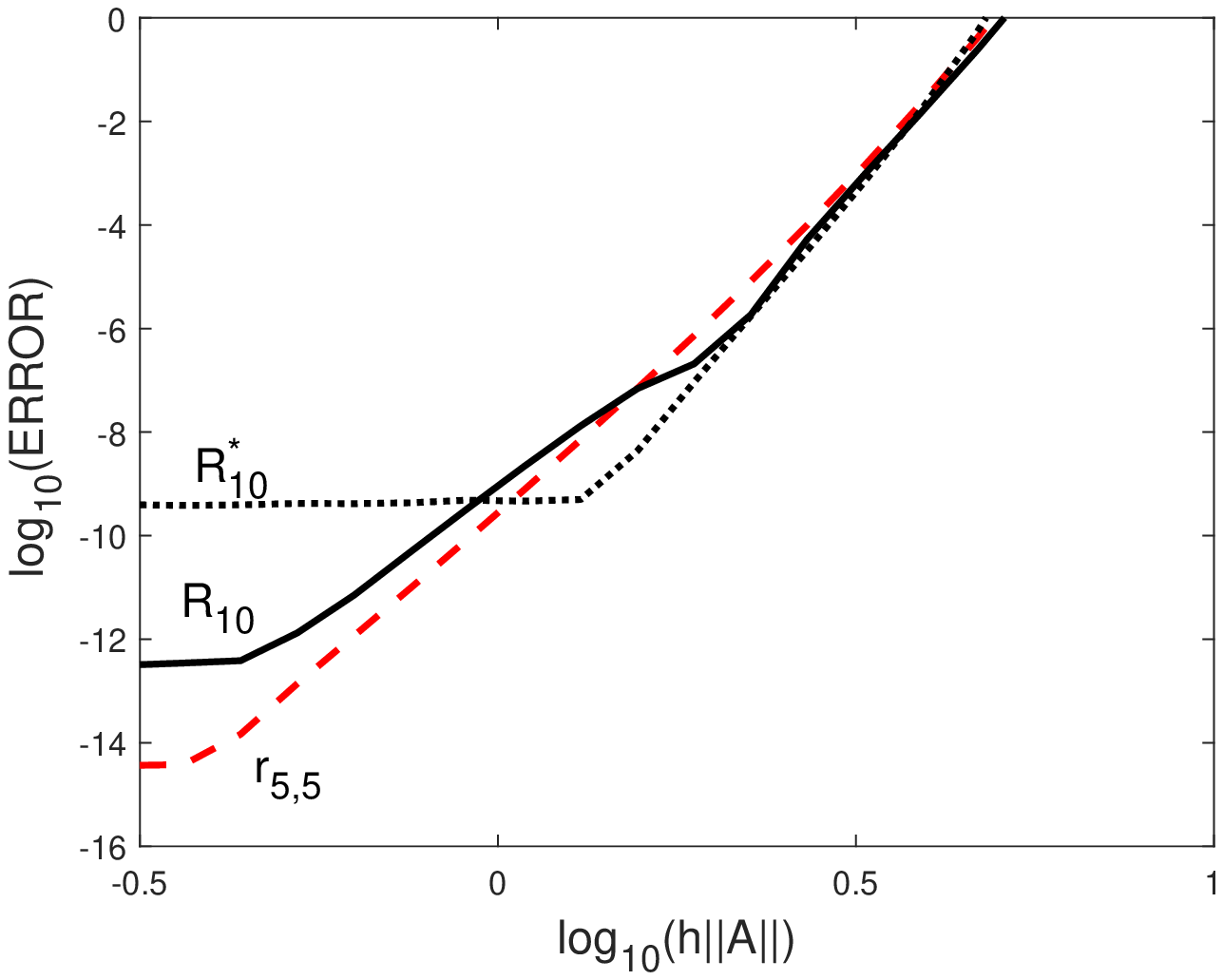}
      \caption{Left: Two-norm error to compute $\e^{hA}$ with $A$ a random matrix of dimension $d=100$ versus $h\|A\|$ in double logarithmic scale. Dashed lines correspond to the results obtained with the 10th-order diagonal Pad\'e method, $r_{5,5}$, doted lines correspond to the 10th-order rational method $R^*_[10]$ and solid lines correspond to the 10th-order rational method $R_[10]$ (solid lines). Right: The same with the matrix \eqref{eq.Ex2b} of dimension $d=100$.}
    \label{Fig.3}
  \end{center}
\end{figure}

We have repeated the numerical experiments using the 10th-order methods.  Figure~\ref{Fig.3} shows the results obtained in a double logarithmic scale:
$r_{5,5}$ (dashed lines), $R^*_{10}$ (dotted lines) and $R_{10}$ (solid lines). We observe the high accuracy of the method $R^*_{10}$, but with large round off errors. The scheme $R_{10}$ shows slightly less accurate results but the round off errors are reduced about three orders of magnitude. Notice that the methods $R^*_{10}$ and $R_{10}$ can be up to three times faster than $r_{5,5}$ and this is not shown in the figure.

As an illustrative example of the action of a function of a large and sparse matrix on a vector we consider the action of the exponential on a vector for the the tridiagonal matrix, 
\[
 A=trid\{-1,2,-1\}, 
\]
of dimension $d=1000$.

We measured the two-norm error of the action of this exponential on a vector,  $\e^{hA}{\bf v}$, for different values of $h$ and where ${\bf v}$ is a unitary vector with random components, i.e. ${\bf w}=\texttt{randn}(d,1), \ {\bf v}={\bf w}/\|{\bf w}\|$ for the 10th-order method, $R_{10}$ and we compare with the Taylor expansion truncated to order 10. 
As previously mentioned, each vector-matrix product involves $5d$ flops while each system requires only 8d flops, so the Taylor method requires 50d flops per step. To reduce round off errors (by nearly two orders of magnitude) one can use the trick given in \eqref{eq.ParallRoundOff} which requires one extra product, i.e. $13d$ flops per step  and processor (approximately 4 times faster). The scheme reads

\begin{equation}\label{eq.Frac10R}
  r_{10}^{(11)}(x)=a_0+d_1x+d_2x^2+\sum_{i=1}^{10} b_i\frac{c_ix}{1-c_ix},
\end{equation}
where we have considered that $d_0+b_1+b_2+\ldots+b_{10}=a_0=1$ so each processor that evaluates each fraction has to compute 
\[
 (I-c_iA){\bf v}_i=c_iA{\bf v}
\]
which requires one product, $A{\bf v}$, and solving a tridiagonal system ($13d$ flops). This method corresponds to $R_{10}$ but we denote it as  $R_{10}^{\prime}$ when the method is computed following this sequence.

\begin{figure}[h!]
  \begin{center}
       \includegraphics[scale=0.6]{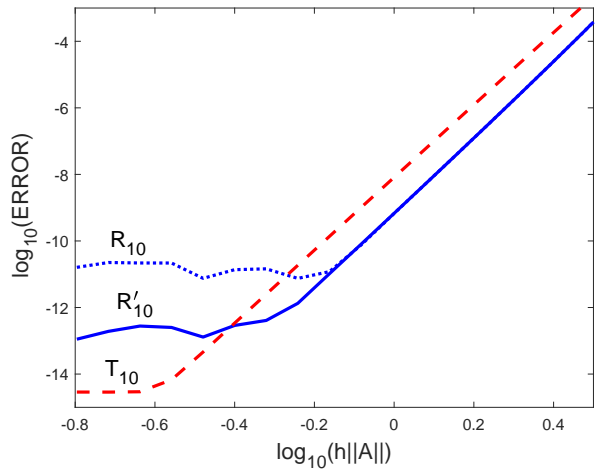}
      \caption{Two-norm error to compute $\e^A{\bf v}$ with $A$ a tridiagonal matrix of dimension $d=1000$ and ${\bf v}$ is a unitary vector with random components versus $h\|A\|$ in double logarithmic scale. Dashed line correspond to the 10th-order Taylor method $T_{10}$, the dotted line to the method $R_{10}$ and the solid line  to the method $R_{10}^{\prime}$.}
    \label{Fig.4}
  \end{center}
\end{figure}

Figure~\ref{Fig.4} shows the results obtained in a double logarithmic scale: dashed line correspond to the 10th-order Taylor method, $T_{10}$, dotted line to the method $R_{10}$ and the solid line $R_{10}^{\prime}$. Notice that  $R_{10}$ and $R_{10}^{\prime}$ are the same method but computed differently which affects both to the round off error as well as to the computational cost. We also remark that $R_{10}$  can be up to more than five times faster than $T_{10}$ and $R_{10}^{\prime}$ up to about four times faster.

\subsection{Exponential of tridiagonal (and banded) matrices acting on vectors}

Let us now consider the particular case (but of great practical interest) of computing the exponential of a tridiagonal matrix acting on a vector, $\e^A{\bf v}$, with single precision accuracy and where we assume that, say $\|A\|\leq 1.5$, as it is usually the case when, for example, exponential integrators are used to solve differential equations. The same algorithms with minor changes in the computational cost of the methods apply to the computation of the exponential of pentadiagonal or banded matrices acting on vectors.

One of the most used schemes to compute the action of the  matrix exponential is proposed in \cite{almohy11ctaf}. The algorithm works as follows: To compute $\e^A{\bf v}$, a set of Taylor polynomials of degree $m=5k, \ k=1,2,\ldots,11$ are considered. Given a tolerance (single or double precision), an estimation to the norm $\|A\|$ (frequently $\|A\|_1$ is used), and from the error analysis, it is chosen the lowest degree Taylor polynomial among the list such that the desired accuracy is guaranteed. In \cite{almohy11ctaf}, the choice of the method is done from the results obtained by the backward error analysis. Very similar results are obtained with the forward error analysis, which is simpler and for the convenience of the reader we will consider it.  In Table~\ref{tab.thetaTay} we collect the relevant values for our problem, i.e. if $\|A\|\leq 0.186$ the algorithm chooses the Taylor approximation of order five, if $0.186 <\|A\|\leq 1.073$ the tenth order method is used, etc. which guarantee that $\|\e^A{\bf v}-T_m(A){\bf v}\|<u\leq 2^{-24}$ where $T_m(A)$ is the Taylor polynomial approximation of degree $m$.

\begin{table}\centering\footnotesize
\caption{\label{tab.thetaTay} Values of $\theta_{m}$ obtained from the forward error analysis for the Taylor approximation $t_m$ of order $m$ for single precision. {The values for the backward error analysis used in the algorithm proposed in \cite{almohy11ctaf} are quite close but slightly shifted, but this also occurs for the rational methods.}
}
\newcolumntype{H}{@{}>{\lrbox0}l<{\endlrbox}}
\newcolumntype{D}{>{$}r<{$}}
\begin{tabular}{DDDDHD}
\toprule
	 m:  & 5					&  			10
	&		  15 
	\\ \midrule	 
	u \leq2^{-24}&
	0.186&
	1.073&
	2.382
	\\
\bottomrule
\end{tabular}
\end{table}

We build a similar scheme in order to analyse the interest of the new algorithms versus the State of the Art algorithm. We then need a 5th-order approximation that, without an exhaustive analysis, we take as the following one:
\[
  r_5^{(5)}(x)=b_1+\sum_{i=2}^6 b_i\frac{1}{1-c_ix}.
\]
If we take $c_1=0$, $c_i=1/(i+1), \ i=2,\ldots,6$ then the system \eqref{eq:LinearEqsbi} has the solution
\begin{equation} \label{eq.Frac5bi0}
  b_1=-\frac{43}{12}, \quad
  b_2=\frac{81}{32}, \quad
  b_3=-\frac{704}{9}, \quad
  b_4=\frac{23125}{48}, \quad
  b_5=-810, \quad
  b_6=\frac{117649}{288}.
\end{equation}
The forward error analysis tells us that the method provides an error below the tolerance,$ u \leq2^{-24}$, for $\|A\|\leq\theta_5=0.298$, as indicated in  Table~\ref{tab.thetaParallel}, being significantly greater than with the Taylor method.
The computational cost to evaluate $r_5^{(5)}(A){\bf v}$ corresponds to solve five tridiagonal systems. If we consider that each system is solved with $8d$ flows, this will correspond to the cost of 8/5 times the cost of the products $A{\bf v}$. If the method is computed sequentially, we observe that the cost would be equivalent to 8 matrix-vector products, and these numbers correspond to the interval $(\frac{8}{5},8)$ shown in  Table~\ref{tab.thetaParallel} as a measure of the cost in comparison with the Taylor method. For the 10th-order method we take the scheme \eqref{eq.Frac10}-\eqref{eq.Frac10biNO} (the scheme with coefficients given in \eqref{eq.Frac10ci}-\eqref{eq.Frac10bi} have smaller round off errors but on the other hand  it is slightly less accurate and and it has a smaller value of $\theta_{10}$). The cost and value of  $\theta_{10}$ for this method are collected in Table~\ref{tab.thetaParallel}.

\begin{table}\centering\footnotesize
\caption{\label{tab.thetaParallel} Values of $\theta_{m}$ obtained from the forward error analysis for the fractional approximations $r^{(5)}_5$ and $r^{(11)}_{10}$ of order $5$ and 10, respectively, for single precision. The cost $m^*$ is measured as an interval in terms of the cost of the product $A{\bf v}$ to makes easier the comparison with the Taylor methods. The lower limit of the interval corresponds to the ideal case where all calculations are carried in parallel and there is no cost to communicate between processors, and the upper limit corresponds to the same scheme fully computed as a serial method.}
\newcolumntype{H}{@{}>{\lrbox0}l<{\endlrbox}}
\newcolumntype{D}{>{$}r<{$}}
\begin{tabular}{DDDDHD}
\toprule
	     m^*:     & (\frac{8}{5},8)	& 	(2,18) 	\\ \midrule	 
	u \leq2^{-24} & 	0.298         & 	1.734
	\\
\bottomrule
\end{tabular}
\end{table}

\begin{figure}[h!]
  \begin{center}
       \includegraphics[scale=0.6]{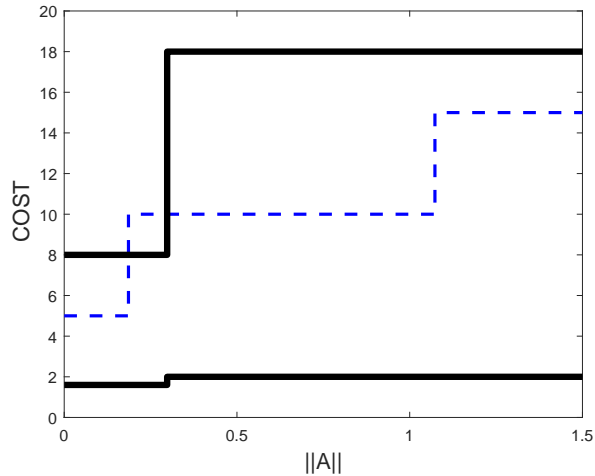}
      \caption{Computational cost to evaluate $\e^A{\bf v}$ with $A$ a tridiagonal matrix (measured as times the cost of one matrix-vector product) versus $\|A\|$ for the algorithm using Taylor approximations of degree 5, 10 and 15 as given in \cite{almohy11ctaf} (dashed line) and the new rational approximations $R_5$ and $R_{10}$ of order $5$ and 10, respectively. The upper solid line corresponds to the worst case in which the cost is measured as if the methods are computed sequentially while the lower one corresponds to the ideal case in which they are computed in parallel with no cost in the communication between processors.}
    \label{Fig.5}
  \end{center}
\end{figure}

The results from Tables~\ref{tab.thetaTay} and \ref{tab.thetaParallel} are illustrated in Figure~\ref{Fig.5} where we plot the computational cost of each algorithm, measured in terms of the cost of a matrix-vector product, for different values of $\|A\|$ in the interval of interest. For the fractional methods (solid lines) the cost is measured both considering they are computed as a sequential scheme as well as if they were computed in an ideal parallel computer (the cost for one single processor). The Taylor methods correspond to the dashed line. We observe that the new methods are competitive even in the worst scenario and much better once some cost is saved due to the parallel programing. The excellent results obtained motivated us to make deeper analysis for this problem in order to find an optimised set of methods for different orders of accuracy similarly to the actual sequential scheme and to test them in multiprocessors workstations. Thsi work will be carried out in the future.

We conclude this section with some remarks:
\begin{itemize}
	\item If one is interested to compute the exponential $\e^A$ acting on several vectors (or the same exponential is applied on several steps for a given integrator) then the $LU$ factorization requires to be done only once making the new algorithms cheaper and competitive even when used as sequential algorithms.
	\item If the matrix $A$ is pentadiagonal and we take into account that the vector-matrix product, $A{\bf v}$, can be done with $9d$ flops and that the pentadiagonal system $(I-c_iA){\bf x}={\bf v}$ can be solved with only $15d$ flops using an $LU$ factorization then, the same conclusions remain approximately valid for this problem because the relative cost is very similar to the case of tridiagonal matrices.
	\item If the matrix $A$ is banded and the linear system  $(I-c_iA){\bf x}={\bf v}$ can be solved accurately with few iterations using, for example, a conjugate gradient method with an incomplete $LU$ factorization as a preconditioner, then the previous results provide a good picture of the performance of the new methods for banded matrices.
	\item The computational cost from the linear combination of vectors should also to be considered in the simple case of the exponential of tridiagonal matrices, making all methods slightly more costly. We have not included this extra cost to provide some results which can also be applied to pentadiagonal or banded matrices where this extra cost is marginal.
\end{itemize}

\section{Conclusions}\label{sec:conclusions}

In this work we have presented a new procedure to compute functions of matrices as well as their action on vectors designed for parallel programming that can be significantly more efficient than existing sequential methods. Given a dense matrix, $A$, the computation of $(I-c_iA)^{-1}$, for sufficiently small constant $c_i$, can be evaluated at the cost of $4/3$ times the cost of a matrix-matrix product and can formally be written as a series expansion that contains all powers of $A$. Then, a proper linear combination of this matrix evaluated in parallel for different values of $c_i$ can allow to approximate any function inside the radius of convergence. If the computation can be carried in parallel with a reduced cost in the communication between processors then the new methods can be up to several times faster than conventional serial algorithms.
For large dimensional problems in which one is only interested in the matrix function acting on a vector, the performance of the new methods depend on the existence of a fast algorithm to solve the linear system of equations,  $(I-c_iA){\bf v}_i={\bf v}$. We have illustrated with some examples that to construct new methods as well as to carry an error analysis is quite simple. The preliminary results are very much promising and it deserves to be further investigated to get optimal methods for different classes of problems and number of processors available in order to get accurate and stable solutions with small round off errors.

\section*{Acknowledgements}
This work was supported by Ministerio de Ciencia e Innovaci\'on (Spain) through 
project PID2019-104927GB-C21 (AEI/FEDER, UE). 



\begin{thebibliography}{00}




\bibitem{almohy09ans} 
{\sc A.~H. Al-Mohy and N.~J. Higham},
 {A new Scaling and Squaring Algorithm for the Matrix Exponential},
 SIAM J. Matrix Anal. Appl., {31}, (2009), pp. 970--989.


\bibitem{almohy11ctaf}
 {\sc A.H. Al-Mohy and N.J. Higham},
Computing the action of the matrix exponential, with an application to exponential integrators,
SIAM J. Sci. Comput. 33 (2011), 488-511.

\bibitem{almohy15naf}
 {\sc A.H. Al-Mohy, N.J. Higham and S.D. Relton},
New algorithms for computing the matrix sine and cosine separately or simultaneously, 
SIAM J. Sci. Comput.  37 (2015)  A456 - A487.


\bibitem{almohy22apa} 
{\sc A.~H. Al-Mohy and N.~J. Higham, and X. Liu},
 {Arbitrary Precision Algorithms for Computing the Matrix Cosine and its Fréchet Derivative},
 SIAM J. Matrix Anal. Appl., {43}, (2022), pp. 233--256.

\bibitem{arioli96tpm} 
{\sc M. Arioli, B. Codenotti, and C. Fassino},
The Pad\'e method for computing the matrix exponential. 
{Lin. Alg. Applic.} {240} (1996), 111--130.


\bibitem{auckenthaler10mea}
{\sc T. Auckenthaler, M. Bader, T. Huckle, A. Sp\"orl, and K. Waldherr},
{Matrix exponentials and parallel prefix computation in a quantum control problem},
Parallel Computing 36 (2010), pp. 359--369.

\bibitem{bader17aia}
{\sc P.~Bader, S.~Blanes, and F.~Casas}, 
{An improved algorithm to compute the exponential of a matrix}.
\newblock arXiv:1710.10989 [math.NA], 2017.

\bibitem{bader19ctm}
{\sc P.~Bader, S.~Blanes, and F.~Casas}, 
{Computing the matrix exponential with an optimized {T}aylor polynomial 
approximation},
  Mathematics, 7 (2019), p.~1174.


\bibitem{bader22aea}
{\sc P.~Bader, S.~Blanes, F.~Casas, and M.~Seydao\u{g}lu}, 
An efficient algorithm to compute the exponential of skew-Hermitian
matrices for the time integration of the Schrödinger equation,
Math. Comput. Sim. 194 (2022), pp. 383--400.


\bibitem{biamonte17qaf} 
{\sc J. Biamonte, P. Wittek, N. Pancotti, P. Rebentrost, N. Wiebe, S. Lloyd},
 {Quantum Machine Learning},
 Nature, {549}, (2017), pp. 195--202.


\bibitem{blanes21npi}
{\sc S.~Blanes}, 
Novel parallel in time integrators for ODEs, 
Appl. Math. Lett., 122 (2021) 107542.





\bibitem{blanes09tme} 
{\sc S.~Blanes, F.~Casas, J.~A.~Oteo, J.~Ros},
{The Magnus expansion and some of its applications},
 {Phys. Rep.}, {470}, (2009), pp. 151--238.



\bibitem{blanes17hoc}
{\sc S.~Blanes, F.~Casas, and M.~Thalhammer}, 
{ High-order commutator-free quasi-magnus exponential integrators for 
non-autonomous linear evolution equations}, 
Comput. Phys. Comm., 220 (2017), pp.~243--262.


\bibitem{calvetti95ipf}
{\sc D. Calvetti, E. Gallopoulos, and L. Reichel},
Incomplete partial fractions for parallel evaluation of rational matrix functions,
J. Comput. Appl. Math., 59 (1995) 349-380.

\bibitem{chehab16pmf}
{\sc J.-P. Chehab and M. Petcu},
Parallel matrix function evaluation via initial value ODE modeling
Comput. Math. Appl., 72 (2016), pp. 76--91.


\bibitem{gallopoulos16pim}
{\sc E. Gallopoulos, B. Philippe, and A.H. Sameh},
Parallelism in matrix computations, 
Springer, 2016.


\bibitem{gallopoulos89otp}
{\sc E. Gallopoulos and Y. Saad},
On the parallel solution of parabolic equations, 
in Proc. 1989 ACM Int’l. Conference on Supercomputing, Herakleion, Greece, June 1989, pp. 17–28. 

\bibitem{gallopoulos92eso}
{\sc E. Gallopoulos and Y. Saad},
Efficient Solution of Parabolic Equations by Krylov Approximation Methods,
SIAM J. Sci. Statist. Comput., 13 (1992), pp. 1236--1264.

\bibitem{gander13pap}
{\sc M.J. Gander and S. G\"uttel},
PARAEXP: A Parallel integrator for linear initial-value problems,
SIAM J. Sci. Comput., 35 (2013), pp. C123--C142.


\bibitem{golub93mc} 
{\sc G.~H. Golub and C.~F. Van Loan},
 {Matrix Computations},
 John Hopkins Univ. Press, Baltimore, Maryland, 1993.


\bibitem{harrow09qaf} 
{\sc A.W. Harrow, A. Hassidim, and S. Lloyd},
 {Quantum Algorithm for Linear Systems of Equations},
 Phys. Rev. Lett., {103}, (2009), 150502.




\bibitem{higham05tsa} 
{\sc N.~J. Higham},
 {The Scaling and Squaring Method for the Matrix Exponential},
 SIAM  J. Matrix Anal. Appl., {26}, (2005), pp. 1179--1193.

\bibitem{higham09tsa} 
{\sc N.~J. Higham},
 {The Scaling and Squaring Method for the Matrix Exponential Revisited},
 SIAM Review, {51}, (2009), pp. 747--764.



\bibitem{higham08fom} 
{\sc N.~J. Higham}, 
{Functions of Matrices: Theory and Computation}, 
Society for Industrial and Applied Mathematics, Philadelphia, PA, USA (2008).



\bibitem{higham10cma} 
{\sc N.~J. Higham and A.~H. Al-Mohy},
 {Computing Matrix Functions},
 Acta Numerica, {51}, (2010), pp. 159--208.


\bibitem{hochbruck97oks}
{\sc M. Hochbruck, Ch. Lubich},
On Krylov subspace approximations to the matrix exponential operator,
SIAM J. Numer. Anal. 34 (1997), 1911-1925

\bibitem{hochbruck10ei}
{\sc M. Hochbruck and A. Ostermann}, 
{Exponential integrators}, 
Acta Numerica, 19 (2010), pp. 209--286.



\bibitem{iserles00lgm} 
{\sc A. Iserles, H.Z. Munthe-Kaas, S.P. N{\o}rsett, and A. Zanna},
Lie group methods. 
{Acta Numerica}  {\em 9} (2000), 215--365.


\bibitem{jarlebring21cgf} 
{\sc E. Jarlebring, M. Fasi, and E. Ringh}, 
Computational graphs for matrix functions,
arXiv preprint arXiv:2107.12198, (2021).



\bibitem{moler03ndw} 
{\sc C.~B. Moler and C.~F. Van Loan},
{Nineteen dubious ways to compute the exponential of a matrix, twenty-five years later},
 SIAM Review {45} (2003), pp. {3--49}.


\bibitem{najfeld95dot} 
{\sc I. Najfeld and T.F. Havel},
Derivatives of the matrix exponential and their computation.
{Adv. Appl. Math.}  {\em 16} (1995), 321--375.


\bibitem{niesen12aak} 
{\sc J. Niesen and W. M. Wright},
Algorithm 919: A Krylov subspace algorithm for evaluating the $\varphi$-functions appearing in exponential integrators, ACM Trans. Math. Software, 38 (2012), Art. 22, 19 pp.




\bibitem{press97nr} 
{\sc  W.H. Press, B.P. Flannery, S.A. Teukolsky, and W.T. Vetterling},
Numerical Recipes in Fortran 90: The Art of Parallel Scientific Computing,
Volume 2 of Fortran Numerical Recipes, Second Ed.,
Cambridge Univ. Press, New York 1997.



\bibitem{sastre18eeo}
{\sc  J.~Sastre}, 
{Efficient evaluation of matrix polynomials}, 
Lin. Alg. Applic., 539 (2018), pp.~229--250.

\bibitem{sastre19btc}
{\sc J.~Sastre, J.~Ib\'a{\~n}ez, and E.~Defez}, 
{Boosting the computation of the matrix exponential}, 
Appl. Math. Comput., 340 (2019), pp.~206--220.


\bibitem{sastre14aae} 
{\sc J. Sastre, J. Ib\'a\~nez, P. Ruiz, and E. Defez},
{Accurate and efficient matrix exponential computation},
Int. J. Comput. Math., 91, (2014), pp. 97--112.

\bibitem{seydaoglu21ctm}
{\sc M.~Seydao\u{g}lu, P.~Bader, S.~Blanes, and F.~Casas}, 
{Computing the matrix sine and cosine simultaneously with a reduced number of products},
  Appl. Numer. Math., 163 (2021), pp.~96--107.



\bibitem{sidje98eas} 
{\sc R.~B.  Sidje},
 {Expokit: a software package for computing matrix exponentials},
 {ACM Trans. Math. Software} {24} (1998), pp. 130--156.


\bibitem{trefethen97nla} 
{\sc L.N. Trefethen and D. Bau III},
 {Numerical linear algebra},
 SIAM, Philadelphia, 1997.



\bibitem{wu21hpc} 
{\sc F. Wu, K. Zhang, L. Zhu, and  J. Hu},
 {High-performance computation of the exponential of a large sparse matrix},
 SIAM  J. Matrix Anal. Appl., {42}, (2021), pp. 1636--1655.



%

	




\end{thebibliography}
\end{document}